\documentclass[12pt,twoside]{article}

\usepackage[T1]{fontenc}
\usepackage[cp1250]{inputenc}

\usepackage[arrow, matrix, curve]{xy}
\usepackage{amsmath,amsthm,amssymb}
\usepackage{amscd,graphics}
\usepackage[dvips]{graphicx}
\usepackage{makeidx}
\newenvironment{Proof}{\textbf{Proof.}}{$\qquad \blacksquare$\par}
\newenvironment{Proof of}[1]{\textbf{Proof #1.}}{$\qquad \blacksquare$\par}

\DeclareMathOperator{\Aut}{Aut}

\DeclareMathOperator{\hull}{hull}
\newcommand{\I}{\mathcal I}

\newcommand{\LL}{\mathcal L}

\newcommand{\TT}{\mathcal T}

\newcommand{\JJ}{\mathcal J}
\newcommand{\KK}{\mathcal K}

\newcommand{\OO}{\mathcal{O}}

\newcommand{\tdelta}{\widetilde \delta}

\newcommand{\FF}{\mathcal F}
\newcommand{\A}{\mathcal A}

\newcommand{\Z}{\mathbb Z}
\newcommand{\N}{\mathbb N}
\newcommand{\T}{\mathbb T}

\newtheorem{thm}{Theorem}[section]

\newtheorem{prop}[thm]{Proposition}
\newtheorem{cor}[thm]{Corollary}
\theoremstyle{definition}

\newtheorem{defn}[thm]{Definition}

\newtheorem{rem}[thm]{Remark}

 \makeindex

\begin{document}

 \thispagestyle{empty}

 \begin{center}
{\bfseries \large \textsc{Relative Cuntz-Pimsner Algebras, \\ Partial Isometric Crossed Products\\
and Reduction of Relations}}

\bigskip
B. K.  Kwa\'sniewski, \ \ A.V. Lebedev  

\end{center}

\begin{abstract} The article discusses the interrelation between relative 
Cuntz-Pimsner algebras and  partial isometric crossed products, and presents a procedure that reduces any given  Hilbert bimodule to  the  "smallest" Hilbert bimodule yielding the  same relative Cuntz-Pimsner algebra as the initial one. In the context of crossed products this reduction procedure corresponds to reduction of $C^*$-dynamical systems.

\end{abstract}

\medbreak


\medbreak
\vspace{5mm}

\tableofcontents

\section*{Introduction}

The Cuntz-Pimsner algebras and relative Cuntz-Pimsner algebras arise in a natural way
as certain generalizations of Cuntz-Krieger algebras \cite{p} so also in the graph algebras theory
(see, for example, the corresponding discussion in \cite{ms}      and \cite{fmr}). By their origin these  algebras are also related 
to crossed products (recall again \cite{p}) which are known to be among the most important structures in 
$C^*$-algebra theory, and until the present day  various crossed product constructions, 
especially those associated with endomorphisms, appear almost continuously (see, for example, the discussion in \cite{Ant-Bakht-Leb} and 
\cite{kwa-leb}). 
If one starts with an arbitrary $C^*$-algebra and its endomorphism then 
the crossed product construction presented in \cite{kwa-leb} can be naturally considered as the most general one, cf. Table 1 
of the present  article. On the appearance of  \cite{kwa-leb} B. Solel noted to the authors that the crossed product constructed 
in \cite{kwa-leb} can  also be modeled as a certain relative Cuntz-Pimsner algebra (Proposition \ref{universality proposition} 
of the present article describes in essence the main idea of B. Solel's remark). Thus we naturally arrive at the discussion of the interrelations between  relative Cuntz-Pimsner algebras and crossed products and this is the theme of the article.
\\
During the discussion  we  take the opportunity  to show that from the point of view of relative Cuntz-Pimsner algebras $\OO(E,J)$ it is enough to consider the case when $E$ is naturally embedded into $\OO(E,J)$ as in other case one may pass to  a smaller \emph{reduced Hilbert bimodule} possessing that property. In particular this allows us to prove that if $E$ is a Hilbert bimodule associated with a $C^*$-dynamical system  then every relative Cuntz-Pimsner algebra $\OO(E,J)$ arises as the crossed product considered  in  \cite{kwa-leb} but  applied to a \emph{reduced $C^*$-dynamical system}.

In the first Section  we recall the construction of relative Cuntz-Pimsner algebras introduced by 
P. S. Muhly and B. Solel in \cite{ms}. The second Section contains the description of the 
reduction procedure of Hilbert bimodules that leads to the "smallest" one 
 giving the same relative Cuntz-Pimsner algebra  as initial one. In the final third Section we recall  
the definition of the crossed products by endomorphisms introduced in \cite{kwa-leb}, 
discuss interrelations between various crossed products presented in Table 1, establish the isomorphism between crossed products and 
relative Cuntz-Pimsner algebras associated to $C^*$-dynamical system Hilbert bimodule (Proposition \ref{universality proposition})
and describe the corresponding reduction procedure for crossed products.

%Thus we see that if $J \cap \ker \delta =\{0\}$, then  the algebra  $\OO(J,E)$ coincides with the crossed product associated with $J$ investigated in \cite{kwa-leb}. Moreover, $\OO(A,E)$ is the isometric crossed product from \cite{Adji_Laca_Nilsen_Raeburn} and $\OO(\{0\},E)$ is the partial-isometric crossed product considered in \cite{Lin-Rae}.    Before we discuss the relationship between  $\OO(J,E)$ and different versions of crossed products in more detail 
\par
\textbf{Conventions.} For the simplicity we shall assume that all the algebras and their representations are  unital.  Throughout the paper $A$ stands for  a unital   $C^*$-algebra and $E$ denotes  a  Hilbert bimodule over $A$, i.e. $E$ is a right Hilbert $A$-module  with the left action  given by a homomorphism $\phi:A\to\LL(E)$  where  $\LL(E)$ is the $C^*$-algebra of adjointable operators  on $E$.  For $x,y \in E$,  we denote by ${\Theta}_{x,y}\in \LL(E)$ the "one-dimensional operator": ${\Theta}_{x,y}(z)=x \cdot \langle y,z\rangle_A$, and  $\KK(E)$ denotes the $C^*$-subalgebra of $\LL(E)$ generated by the operators ${\Theta}_{x,y}$, $x,y\in E$.  
\section{Toeplitz representations and associated algebras}

Suppose that $A$ is a unital $C^*$-algebra and $E$ is a Hilbert bimodule over $A$, where the left action $a\cdot x$ is given by a homomorphism $\phi:A\to\LL(E)$, so that $a\cdot x=\phi(a)x$. 
A {\em Toeplitz representation} $(\psi,\pi)$ of the Hilbert bimodule
$E$
 in a unital
$C^*$-algebra
$B$ consists of a linear map $\psi:E\to B$
and a unital homomorphism $\pi:A\to B$ such that
\begin{align}
\psi(x\cdot a) &= \psi(x)\pi(a),\label{eq:rep1} \\
\psi(x)^*\psi(y)&= \pi(\langle x,y\rangle_A),\ \mbox{ and}\label{eq:rep2} \\
\psi(a\cdot x) &= \pi(a)\psi(x).\label{eq:rep3}
\end{align}
for $x,y\in E$ and $a\in A$. We  recall \cite[Remark 1.1]{fr}.
 
\begin{rem}\label{remark:psi}  Let us note that condition \eqref{eq:rep2}
itself already implies that $\psi$ is linear.
It also implies that $\psi$ is bounded: for $x\in E$ we have
$$
\|\psi(x)\|^2 =\|\psi(x)^*\psi(x)\| =\|\pi(\langle x,x\rangle_A)\|
\le \|\langle x,x\rangle_A\| = \| x\|^2.
$$
If $\pi$ is injective, then we have equality throughout, and $\psi$
is isometric.
\end{rem}
 Given a Toeplitz representation $(\psi,\pi)$ , \cite[Proposition~1.6]{fr} says there is a homomorphism ${(\psi,\pi)}^{(1)}:\KK(E)\to B$ which satisfies
$$
{(\psi,\pi)}^{(1)}({\Theta}_{x,y})=\psi(x){\psi(y)}^*\text{ for } x,y\in E,
$$
and
$$
{(\psi,\pi)}^{(1)}(T)\psi(x)=\psi(Tx)\text{ for }T\in\KK(E)\text{ and }x\in E.
$$
We define
\[
J(E):={\phi}^{-1}(\KK(E)),
\]
which is a closed two-sided ideal in $\A$. Let $J$ be an ideal contained in $J(E)$. We say that a Toeplitz representation $(\psi,\pi)$ of $E$ is {\em coisometric on} $J$ if
\begin{equation*}
{(\psi,\pi)}^{(1)}(\phi(a))=\pi(a)\quad\text{for all \ } a\in J.
\end{equation*}

\subsection{The Fock representation }\label{The Toeplitz C*-algebra of a Hilbert bimodule}

Given a Hilbert bimodule $E$ over $A$, for $n\ge 1$, the $n$-fold internal tensor product
$E^{\otimes n} := E\otimes_A\dotsm\otimes_A E$ is naturally a right
Hilbert $A$-module, and $A$ acts on the left by
\[
\phi^{(n)}(a)(x_1\otimes_A\dotsm\otimes_A x_n) :=
(a\cdot x_1)\otimes_A\dotsm\otimes_A x_n;
\]
 For $n=0$, we take $E^{\otimes 0}$
to be the Hilbert module
$A$ with left action
$\phi^{(0)}(a) b: = ab$. Then the Hilbert-module direct sum 
$$
\FF(E) :=
\bigoplus_{n=0}^\infty E^{\otimes n}
$$ carries a diagonal left action $\phi_\infty$
of $A$ in which $\phi_\infty(a)(x):=\phi^{(n)}(a)x$ where  $x\in E^{\otimes n}$.
The Hilbert bimodule $\FF(E)$ is called 
%in \cite{ms} 
 the \emph{Fock space} over  the Hilbert bimodule $E$. 
For each  $x\in E$, we 
define a \emph{creation operator\/}
$T(x)$  on $\FF(E)$ by
$$
T(x)y=\begin{cases}
 x\cdot y
    & \text{if $y\in E^{\otimes 0}=A$} \\
  x\otimes_A y
    & \text{if $y\in E^{\otimes n}$ for some $n\geq 1$;} \\
\end{cases}
$$
routine calculations show that $T(x)$ is adjoint to the  \emph{annihilation operator}
$$
T(x)^*z=\begin{cases}
 0
    & \text{if $z\in E^{\otimes 0}=A$} \\
  \langle x, x_1\rangle_A\cdot y
    & \text{if $z=x_1\otimes_A y\in E\otimes_A E^{\otimes n-1}=
E^{\otimes n}$.}
\\
\end{cases}
$$

It is clear that $T:E\to\LL(\FF(E)) $ is an injective linear mapping and since $A$ is a summand of $\FF(E)$, the map $\phi_\infty:A \to \LL(\FF(E))$ is injective as well. Moreover, the pair is  $(T,\phi_\infty)$ is a Toeplitz representation of
$E$.
\begin{defn}
The Toeplitz representation   $(T,\phi_\infty)$  of
$E$ in the $C^*$-algebra
$\LL(\FF(E))$ is called \emph{Fock representation} and the  \emph{Toeplitz} $C^*$-algebra $\TT(E)$ of $E$ is by definition the $C^*$-subalgebra of
$\LL(\FF(E))$ generated by
$T(E)\cup\phi_\infty(A)$, cf. \cite[Definition 2.4]{ms}, \cite[Definition 1.1]{p}.
 \end{defn}

% A homomorphism ${(T,\phi_\infty)}^{(1)}:\KK(E)\to\LL(\FF(E))$ is injective and since $$ {(T,\phi_\infty)}^{(1)}({\Theta}_{x,y})=T(x){T(y)}^*\in \TT(E)\qquad \text{ for } x,y\in E, $$we may view $\KK(E)$ as isometrically contained in $\TT(E)$. \\
\subsection{Relative Cuntz-Pimsner algebras $\OO(E,J)$}\label{Relative Cuntz-Pimsner algebras}

 Let $J$ be an   ideal  in $A$ contained in $J(E)={\phi}^{-1}(\KK(E))$ and let $P_0$ be the projection in $\LL(\FF(E))$ 
 that maps $\FF(E)$ onto the first summand $A$. 
 One can  show \cite[Lemma 2.17]{ms} that $\phi_\infty(J)P_0$  is contained in $\TT(E)$. 
 We shall  write $\JJ(J)$ for the ideal in $\TT(E)$ generated by $\phi_\infty(J)P_0$. 
 
 \begin{defn}
If $E$ is a Hilbert bimodule over the $C^*$-algebra $A$, and if  $J$ is an ideal in  $J(E)={\phi}^{-1}(\KK(E))$ we denote by $\OO(J,E)$ the quotient algebra $\TT(E)/\JJ(J)$ and  call it \emph{relative Cuntz-Pimsner algebra} determined by $J$.
 \end{defn}
 The Cuntz-Pimsner algebra $\OO(J,E)$ is universal with respect to Toeplitz representations that are  coisometric on $J$ in the following sense, see  \cite[Proposition~1.3]{fmr}.
\begin{prop}\label{RCP algebra}
Let $E$ be a Hilbert bimodule over $A$, and let $J$ be an ideal in $J(E)$. Let $q:\TT(E)\to \OO(J,E)$ be the quotient map and put 
$$
k_E=q \circ T \quad\textrm{ and } \quad k_A= q\circ \phi_\infty.
$$
 Then  $(k_E,k_A)$ is  a Toeplitz representation of $E$ which is coisometric on $J$ and satisfies:
\begin{itemize}
\item[(i)] for every Toeplitz representation $(\psi,\pi)$ of $E$ which is coisometric on $J$, there is a homomorphism $\psi{\times}_J\pi$ of $\OO(J,E)$ such that 
$$(\psi{\times}_J\pi)\circ k_E=\psi \quad \textrm{  and } \quad (\psi{\times}_J\pi)\circ k_\A=\pi,$$ 
\item[(ii)] $\OO(J,E)$ is generated as a $C^*$-algebra by $k_E(E)\cup k_A(A)$.
\end{itemize}
The triple $(\OO(J,E),k_E,k_A)$ is unique in the following sense: if $(B,k_E',k_A')$ has similar properties, there is an isomorphism $\theta:\OO(J,E)\to B$ such that $\theta\circ k_E=k_E'$ and $\theta\circ k_\A=k_\A'$. There is a strongly continuous gauge action $\gamma:\T\to\Aut\OO(J,E)$ which satisfies ${\gamma}_z(k_A(a))=k_\A(a)$ and ${\gamma}_z(k_E(x))=zk_E(x)$ for $a\in A,x\in E$.
\end{prop}

The algebra $\OO(\{0\},E)$ is the Toeplitz algebra $\TT(E)$ and $\OO(J(E),E)$ is the Cuntz-Pimsner algebra $\OO(E)$ \cite{p}.
The following proposition tells us when $k_A:A\to\OO(J,E)$  is injective (if it is so then $k_E$ is also injective, cf.  Remark \ref{remark:psi}), see  \cite[Proposition~2.21]{ms} and \cite{Brow-Rae}. 

\begin{prop}\label{injectivity of k_A}
Let $E$ be a Hilbert bimodule over $A$ and let $(\OO(J,E),k_A,k_E)$ be a relative Cuntz-Pimsner algebra associated to $E$. Then $k_A$ is injective if and only if 
\begin{equation}\label{orthogonal property}
\ker\phi \cap J =\{0\}.
\end{equation}
\end{prop}
Let us  note that for any ideal $I$  in $A$ the family of ideals $J$ in $A$ such that 
$
I\cap J =\{0\},
$
possess the largest element $I^{\bot}$ (in the sense of partial order given by inclusion).  Namely, we have
$$
I^\bot = \bigcap_{x\notin  \hull(I)}\,  x.
$$
where  $
\hull(I)=\{x\in {\rm Prim} A: \ x\supset I \}
$ is the hull of $I$.
\\
Within this notation, relation \eqref{orthogonal property} is equivalent to the inclusion $
J\subset (\ker\phi)^{\bot}$. The aim of the present paper is to show that  all the relative Cuntz-Pimsner algebras are the algebras $\OO(J,E)$ where   
$J\subset (\ker\phi)^{\bot}$.   Moreover, if the choice of an ideal $J$ is not dictated by any outside demands it seems that 
$$
J=(\ker\phi)^{\bot}\cap J(E)
$$ 
 is the best one to choose.

\section{Reduction of Hilbert bimodules}\label{Reduction of Hilbert bimodule}
We recall certain results from \cite{fmr} concerning invariant ideals and quotient bimodules.
\\
Let $E$ be a bimodule over $A$ and let $I$ be an ideal in $A$. The closed subspace
$$
EI:=\{x \cdot i: x\in E, \,\, i\in I\}
$$
is a right Hilbert $I$-module. Moreover if we let $q^I:A\to A/I$ and $q^{EI}: E \to EI$ be the quotient maps then by \cite[Lemma 2.1]{fmr}
$E/EI$ is a right Hilbert $A/I$-module  with 
\begin{equation}\label{right action}
q^{EI}(x)\cdot q^I(a):=q^{EI}(x\cdot a),\qquad x \in E,\,\, a\in A,
\end{equation}
\begin{equation}\label{scalar product}
\langle q^{EI}(x),q^{EI}(y)\rangle_{A/I}:=q^{I}(\langle x,y\rangle_A).
\end{equation}
In order to define a left action on $E/EI$ we need to impose on the ideal $I$ that 
\begin{equation}\label{E-invariance}
\phi(I)E\subset EI.
\end{equation}
An ideal $I$ in $A$ satisfying \eqref{E-invariance} is called  $E$-\emph{invariant} and if $I$ is $E$-invariant then by \cite[Lemma 2.3]{fmr} there is  a homomorphism $\phi_{A/I}:A/I\to  \LL(E/EI)$ such that 
\begin{equation}\label{left action}
\phi_{A/I}(q^{I}(a))q^{EI}(x)=q^{EI}(\phi(a)x),\qquad x \in E,\,\, a\in A.
\end{equation}
Thus, for any  $E$-invariant ideal $I$ in $\A$ the space $E/EI$ together with  \eqref{right action},  \eqref{scalar product}  \eqref{left action} is a bimodule. We shall call it \emph{quotient bimodule} of $E$.   
\\
We recall the main theorem from \cite{fmr}.
\begin{thm} \label{takie tam aa}\cite[Theorem 3.1]{fmr}
Suppose $E$ is a Hilbert bimodule over $A$, $J$ is an ideal in $J(E)$, and $I$ is an $E$-invariant ideal in $A$. If we denote by  $\I(I)$ the  ideal in $\OO(J,E)$ generated by $k_A(I)$ then the quotient $\OO(J,E)/\I(I)$ is canonically isomorphic to $\OO(q^I(J),E/EI)$.
\end{thm} 
Let us now fix a Hilbert bimodule $E$ over $A$ and an ideal $J$ in $J(E)$. We will now reduce $E$ by taking quotient of it to a certain "smaller" Hilbert bimodule satisfying \eqref{orthogonal property} and yielding  the same relative  Cuntz-Pimsner algebra as $E$ and $J$. 
\\
We define recursively a sequence of ideals in $A$ putting
$$\label{ideals I_n}
J_0=\{a\in J: \phi(a)=0\}=\ker \phi\cap J,
$$
and for $n\geq 0$,
$$
 J_{n+1}=\{a\in J: \phi(a)E\subset E{J_n}\}.
$$
Then one easily sees that $\{J_n\}_{n\in \N}$ is an  increasing family of $E$-invariant ideals in $A$ and hence
$$
J_\infty=\overline{\bigcup_{n\in \N} J_n}
$$ 
is an $E$-invariant ideal in $A$.  We note that 
\begin{equation}\label{I_infty condition}
a \in J\,\, \wedge \,\, \phi(a)E\subset E J_\infty \Longrightarrow  a\in J_\infty 
\end{equation}
and  this implication characterizes  $J_\infty$ in the sense that it is the smallest $E$-invariant ideal in $A$ containing $\ker \phi\cap J$ and satisfying \eqref{I_infty condition}.\\
The first of our results  states that the quotient Hilbert bimodule $E/EJ_{\infty}$ and the quotient $C^*$-algebra  $A/J_\infty$ may be identified with the image of the initial Hilbert bimodule $E$ and $C^*$-algebra $A$ in the relative Cuntz-Pimsner algebra $\OO(E,J)$. 
\begin{thm}\label{reduction thm}
Let  $E$ be a Hilbert bimodule over $A$ and  $J$  an ideal in $J(E)$. Then for $n\in \N\cup\{\infty\}$, we have a canonical isomorphism  
$$
\OO(J,E) 
\cong \OO(q^{J_n}(J),E/EJ_{n}).
$$
Moreover, for $n=\infty$ we have
$$
\ker \phi_{A/J_\infty}\cap q^{J_\infty}(J)=\{0\}
$$
and thus we have the following (again canonical) isomorphisms 
$$
k_A(A)\cong A/J_\infty, \qquad k_E(E)\cong E/EJ_\infty.
$$
\end{thm}
\begin{Proof} In view of Theorem  \ref{takie tam aa} to prove the first part of theorem it is enough to show that  for every ideal $J_n$, $n=0,1,..., \infty$,  we have  $k_A(J_n)=0$.
\\
 It is  clear that $k_A(J_0)=0$. Assume that $k_A(J_n)=0$ and let $a\in J_{n+1}$. Then for every $x\in E$ there exists $y(x)\in E$ and $i(x)\in J_N$ such that $\phi(a)x=$ and thus 
$$
k_A(a)k_E(x)=k_E(\phi(a)x)= k_E(y(x)i(x))=k_E(y(x))k_A(i(x))=0.
$$
Hence $k_A(J_{n+1})=0$.  It follows that $k_A(J_n)=0$ for every $n=0,1,..., \infty$.
\\
To prove that  $
\ker \phi_{A/J_\infty}\cap q^{J_\infty}(J)=\{0\}
$ take $a\in J$ and suppose that $\phi_{A/J_\infty}( q^{J_\infty}(a))=0$.  Then by  \eqref{left action}  we 
see that $\phi(a) E \subset EJ_{\infty}$ and  by  \eqref{I_infty condition} we have $q^{J_\infty}(a)=0$.
 Now it suffices to apply Proposition \ref{injectivity of k_A} and Remark \ref{remark:psi}.
\end{Proof}

 Ideal $J_\infty$ plays the role of a certain  "measure" of the degree of degeneracy of $\OO(J,E)$ since  the bigger $J_\infty$ is   
  the smaller $\OO(J,E)$ is. 
 In particular, $\OO(J,E)=0$ if and only if $J_\infty=
A$,  and if $J\neq A$, then  $\OO(J,E)\neq 0$. Obviously,  $A$ embeds into $\OO(J,E)$ if and only if $J_\infty=0$ which is  equivalent to $E=E/EJ_{\infty}$.

 Theorem \ref{reduction thm} shows in fact that one may always  restrict his interest only to the relative  Cuntz-Pimsner algebras $\OO(E,J)$ determined by ideals such that 
$$
 J \subset (\ker \phi)^\bot,
$$ 
since in any case one may pass to the \emph{reduced Hilbert bimodule} $E/EJ_\infty$.

Note also that  Hilbert bimodules  $E / EJ_{n}$, $n\in \N$ may be considered as  "approximations" of $E / EJ_{\infty}$. In particular, if $J_n=J_{n+1}$, for certain $n\in \N$, then $J_\infty=J_n$.

\section{Crossed products by endomorphisms and their canonical $C^*$-dynamical systems}\label{Crossed products by endomorphisms}
Let $\delta$ be an endomorphism of a unital $C^*$-algebra $A$.  Throughout the paper the pair $(\A,\delta)$ will be called a {\em $C^*$-dynamical system}. We slightly extend a definition from \cite{kwa-leb}.

\begin{defn}\label{kowariant  rep defn}
Let $(A,\delta)$ be a $C^*$-dynamical system. A \emph{covariant representation} of  $(A,\delta)$ in a $C^*$-algebra $B$ is a doublet $(\pi,U)$ consisting of a unital homomorphism 
$\pi:A\to B$ and an operator $U$ satisfying the following relations
\begin{equation}\label{covariance rel1}
U\pi(a)U^* =\pi(\delta(a)),\qquad a \in A,
\end{equation}
\begin{equation}\label{covariance rel2}
  U^*U \in \pi(A)'.
\end{equation}
Every covariant representation $(\pi,U)$ defines an ideal $J$ in $A$ given by  
\begin{equation}\label{covariance rel3}
J=\{ a\in A: U^*U \pi(a)=\pi(a)\}
\end{equation}
If $J$ is an ideal in $A$ and  $(\pi,U)$ is a covariant representation of $(A,\delta)$ satisfying \eqref{covariance rel3} then we say that $(\pi,U)$  is \emph{associated with} $J$.
\end{defn}

Let us note that \eqref{covariance rel1} implies that $U$ is a partial isometry, and \eqref{covariance rel1}
together with \eqref{covariance rel2} imply  that $U$ is power partial isometry. Moreover, see \cite[Theorem 1.6]{kwa-leb} and  remark below,  a covariant representation $(\pi,U)$ such that $\pi$ is faithful exists if and only if it associated with an ideal $J$ having a zero intersection with the kernel of $\delta$. Thus, if $J\cap \ker \delta=\{0\}$ relations \eqref{covariance rel1}, \eqref{covariance rel2}, \eqref{covariance rel3} give rise to a non-degenerate  universal algebra, which was investigated in \cite{kwa-leb}.  To be more precise, a $C^*$-algebra $C^*(\A,\delta,J)$ introduced in  \cite[Definition 4.2]{kwa-leb} is a $C^*$-enveloping algebra of a certain Banach $^*$-algebra however in view of  \cite[Theorem 5.4]{kwa-leb} and   Proposition \ref{star property} we prove below,   this definition is equivalent to the following one.
\begin{defn}\label{crossed product defn}
Let $(A,\delta)$ be a $C^*$-dynamical system and $J$ an ideal in $A$ such that $J\cap \ker \delta=\{0\}$. A \emph{crossed product} $
C^*(A,\delta,J)
$ of $\A$ by $\delta$ associated with $J$ is   a universal $C^*$-algebra generated by the the copy of the algebra $A$ and a partial isometry $u$ subject to relations 
$$
ua u^*=\delta(a),\,\,\,\qquad  u^*ua=a u^*u,\qquad\quad   a \in A, 
$$
$$
J=\{a\in A: u^*u a=a\}.
$$
\end{defn}
  The proof of the next statement is a standard argument, cf. \cite{Ant-Bakht-Leb}.
 \begin{prop}\label{star property}
  Let $(A,\delta)$ be a $C^*$-dynamical system and $J$ such that $J\cap \ker \delta=\{0\}$. The crossed product  $ C^*(A,\delta,J)$ possess the so-called  $(^*)$-property, that is the following inequality holds
  $$
\|\sum_{m=0}^N u^{*m}a_{m,m}u^m \| \leq \|\sum_{m,n=0}^N u^{*m}a_{m,n}u^n \|,  \qquad\qquad (*)
$$
 where  $a_{m,n}\in A$, $n,m=0,1...,N$, and $N\in \N$.
  \end{prop}
\begin{Proof}
Take any faithful representation $\widetilde{\pi}: C^*(A,\delta,J) \to \LL(H)$ of $ C^*(A,\delta,J)$ on a Hilbert space $H$ and "disintegrate"  $\widetilde{\pi}$ to $(\pi,U)$ where $\pi:=\widetilde{\pi}|_A$ and $U:=\widetilde{\pi}(u)$. Then $(\pi,U)$  is  a covariant representation  of $(A,\delta)$ associated with $J$. Consider the space ${\cal H} =l^2 ({\mathbb Z}, H)$ and the
representation  $\nu  : C^*(A,\delta,J) \to \LL({\cal
H})$ given by the formulae
\begin{gather*}
(\nu (a)\xi )_n = \pi (a) (\xi_n), \qquad\text{where}\quad
a\in A, \quad l^2 ({\mathbb Z}, H) \ni \xi = \{ \xi_n  \}_{n\in
\Z}\,;\\[6pt]
(\nu (u)\xi )_n = U (\xi_{n-1}),\qquad
(\nu (u^*)\xi )_n = U^* (\xi_{n+1}).
\end{gather*}
Routine verification shows that $(\nu|_A, \nu
(u))$ is a covariant representation of $(A,\delta)$ associated with $J$ and thus $\nu$ is indeed a representation of   $C^*(A,\delta,J)$.
\\
Now take any $x=\sum_{m,n=0}^N u^{*m}a_{m,n}u^n\in C^*(A,\delta,J)$ where   $a_{m,n}\in A$, $n,m=0,1...,N$, and $N\in \N$. For a given $\varepsilon > 0$ we may chose  a vector
$\eta \in H$ such that $\Vert \eta \Vert =1$ and  
\begin{equation}\label{e}
\Vert \sum_{m=0}^N U^{*m}\pi(a_{m,m})U^m
\eta \Vert > \Vert\sum_{m=0}^N U^{*m}\pi(a_{m,m})U^m\Vert - \varepsilon .
\end{equation}
Set $\xi = \{ \xi_n  \}_{n\in
\Z} \in l^2 ({\mathbb Z}, H)$ by $\xi_n = \delta_{0n}\eta $,
where $\delta_{ij}$ is the Kronecker symbol. We have that $\Vert
\xi \Vert = 1$ and the explicit form of $\nu (x) \xi$  and
\eqref{e} imply
\[
\Vert \nu (x) \xi \Vert \ge \|\nu(\sum_{m,n=0}^N u^{*m}a_{m,n}u^n ) \xi \|=\Vert \sum_{m=0}^N U^{*m}\pi(a_{m,m})U^m
\eta \Vert  
\]
which by \eqref{e} and the arbitrariness of $\varepsilon$ proves the desired
inequality
\[
\Vert x  \Vert \ge \Vert\sum_{m=0}^N U^{*m}\pi(a_{m,m})U^m\Vert = \Vert \widetilde{\pi}(\sum_{m=0}^N u^{*m}a_{m,m} u^m)\Vert = \Vert \sum_{m=0}^N u^{*m}a_{m,m} u^m\Vert.
\]
\end{Proof}

Proposition \ref{star property}  and
\cite[Theorem 5.4]{kwa-leb} imply that  the crossed products introduced in Definition \ref{crossed product defn} and \cite[Definition 4.2]{kwa-leb}  are canonically isomorphic. 

In general, as it was communicated to authors by B. Solel, algebras of this sort can be also modeled out  as 
certain relative Cuntz-Pimsner algebras of P. S.  Muhly and B. Solel \cite{ms}. 

Indeed, let $(A,\delta)$ be a $C^*$-dynamical system  and define the structure of a  Hilbert bimodule   over $A$ on the space 
$$
E:=\delta(1) A
$$
 by
$$
a \cdot x :=\delta(a)x,\quad x\cdot a:= xa,\quad \textrm{ and}\quad \langle x,y\rangle_A:=x^*y.  
$$
Then one easily  checks that $J(E)=A$  and  $\ker\delta =\ker\phi$. We shall say that $E$ is the \emph{$C^*$-dynamical system Hilbert bimodule} of $(\A,\delta)$. The proof of the foregoing proposition  in essence follows the  argument from \cite[Example 1.6]{fmr}.
\begin{prop}\label{universality proposition}
Let $E$ be  a $C^*$-dynamical system Hilbert bimodule of  $(\A,\delta)$ and let $J$ be an ideal in $A$. The relations 
$$
U=\psi(\delta(1))^*, \qquad \psi(x)= U^*\pi(x)
$$
establish a one-to-one correspondence between  Toeplitz representations of $(\psi,\pi)$ of $E$ which are coisometric on $J$, and covariant representations $(\pi,U)$ which are associated with an ideal containing $J$. 
In particular,   
\begin{itemize}
\item[(i)] $\OO(J,E)$ is generated as a $C^*$-algebra by the partial isometry $u=k_E(\delta(1))^*$ and the $C^*$-algebra $k_A(A)$.
\item[(ii)] for every covariant  representation $(\pi,U)$ of $(A,\delta)$ associated with an ideal containing $J$, there is a homomorphism $\pi{\times}_J U$ of $\OO(J,E)$ uniquely determined by 
$$
(\pi{\times}_J U)(u)= U\qquad \textrm{ and  }\qquad (\pi{\times}_J U)\circ k_A=\pi.
$$
\end{itemize}
\end{prop}
\begin{Proof} Let $(\psi,\pi)$ be a Toeplitz representations of  $E$ and let $U:=\psi(\delta(a))$. Then
$$
\pi(\delta(a))=\pi(\langle \delta(1),\delta(a)\rangle_A)=\psi(\delta(1)^*\psi(\delta(a))=U \psi(a \cdot \delta(1))=U\pi(a)U^*.
$$ 
Another computation shows that
$$
U^*\pi(\delta(a))=\psi(\delta(1)\delta(a))=\psi(\delta(a)\delta(1))=\pi(a)\psi(\delta(1))=\pi(a)U^*,
$$
which implies that $U^*U \in \pi(A)'$, cf. \cite[Proposition 2.2]{Leb-Odz}, \cite[Lemma 4.3]{Lin-Rae}. Thus $(\pi,U)$ is a covariant representation of $(A,\delta)$. Moreover, observe that the operator $\phi(a)$ is just $\Theta_{\delta(a),\delta(1)}$ and thus for every $a\in J$ we have
$$
\pi(a)=\pi^{(1)}(\Theta_{\delta(a),\delta(1)})=\psi(\delta(a))\psi(\delta(1))^*=U^*\pi(\delta(a))U
=U^*U\pi(a)U^*U=U^*U\pi(a)
$$
which means that $(\pi,U)$ is associated with   a certain ideal containing $J$.
\\
Conversely let $(\pi,U)$ be a covariant representation  of $(A,\delta)$  associated with   a certain ideal containing $J$ and let $\psi(x):=V^*\pi(x)$. Then one easily checks the conditions \eqref{eq:rep1}, \eqref{eq:rep2} \eqref{eq:rep3},
and to show that $(\pi,\psi)$ is coisometric on $J$ it is enough to reverse the argument we used above. Indeed, for any $a\in J$ we have
$$
\pi(a)=U^*U\pi(a)=U^*U\pi(a)U^*U=U^*\pi(\delta(a))U =\psi(\delta(a))\psi(\delta(1))^*
$$
$$
= \pi^{(1)}(\Theta_{\delta(a),\delta(1)})=\pi^{(1)}(\phi(a)),
$$
and the proof is complete.
\end{Proof}
By the universality of $\OO(E,J)$ and $C^*(A,\delta,J)$ we get the following 
\begin{cor}
Let $E$ be  a $C^*$-dynamical system Hilbert bimodule of  $(\A,\delta)$ and let $J$ be an ideal in $A$ such that $\ker\delta\cap J=\{0\}$. Then algebras $\OO(E,J)$ and $C^*(A,\delta,J)$ are canonically isomorphic.
\end{cor}

In view of the above statements  relative Cuntz-Pimsner algebra $\OO(J,E)$ seems to be a natural candidate for a crossed product of $(\A,\delta)$ associated with an arbitrary ideal $J$.  
In fact, most of the considered  crossed products by endomorphisms coincide with   $\OO(J,E)$ for certain  $J$.   Table 1  presents the corresponding  juxtaposition of the  objects chosen. To see the coincidence in N.3 of Table 1 we refer the reader to \cite[Proposition 2.6]{kwa-ext}. 
We stress once again that the kernel of the  left action in a $C^*$-dynamical system Hilbert bimodule $E$ coincide with $\ker \delta$ 
and hence in view of Proposition  \ref{injectivity of k_A}, algebra $A$ embeds naturally into $\OO(J,E)$ if and only if $J\cap \ker \delta =\{0\}$, 
or equivalently  
\begin{equation}\label{inclusion for an ideal J}
\{0\}\subset J \subset (\ker\delta)^\bot.
\end{equation}
The crossed product N.5 (see Definition	\ref{crossed product defn}) is the most general in the sense that it gives all the remaining ones   for an appropriate choice of $J$ ($J=(\ker\delta)^\bot$ for N.1-5 and $J=\{0\}$ for N.7).  However, in order to get the crossed product N.6 from N.5  one first  needs to 'reduce' the initial $C^*$-dynamical system (this reduction agrees with  the one discussed in Section \ref{Reduction of Hilbert bimodule}, see Corollary \ref{corl for reduced dynamics}).
\\
We  also have to note that two kinds of  crossed products introduced by R. Exel in \cite{exel1} and   \cite{exel2} respectively, also arise as relative Cuntz-Pimsner algebras but  for  differently defined  Hilbert bimodules, see \cite[Example 2.22]{ms} and  \cite{Brow-Rae}. Furthermore these crossed products may be obtained from  the crossed product N.6 of Table 1, see  \cite{Ant-Bakht-Leb}.
\par
Looking at Table 1 one can not help  feeling that among the ideals satisfying \eqref{inclusion for an ideal J} the ideal $J=(\ker\delta)^\bot$ is somewhat privileged. It is completely natural as   $\OO((\ker\delta)^\bot,E)$ should be considered as 'the smallest' relative Cuntz-Pimsner algebra containing all the information about the $C^*$-dynamical system $(\A,\delta)$. Moreover, we shall show that for an arbitrary choice of $J$ the algebra   $\OO(J,E)$ coincides with  Cuntz-Pimsner algebra  $\OO((\ker\delta_J)^\bot,E_J)$ where $E_J$ is a $C^*$-dynamical system Hilbert bimodule of  a canonically constructed pair $(A_J,\delta_J)$.
\begin{table}
\begin{center}
		\begin{tabular}{|c|c|c|c|} \hline
			N. & endomorphism $\delta:A\to A$ & $J \triangleleft A$ & $\OO(J,E)$  
			\\ \hline
			1. & automorphism    & $J=(\ker\delta)^\bot=A$                 & classical unitary 
			\\ &                 &                                        & crossed product 
			\\ \hline
			2. & monomorphism   & $J=(\ker\delta)^\bot=A$                 & isometric crossed product
			\\ &                &                                         &    \cite{Paschke}, \cite{Murphy}   
			\\ \hline
			3. & $\ker\delta$ unital and     & $J=(\ker\delta)^\bot$                 & crossed product using
			\\ &  $\quad \delta(A)$ hereditary in $A\quad $    &                     & complete transfer operator 
			\\ &                                               &                     & \cite{Ant-Bakht-Leb}      
			\\ \hline
		  4. & $\ker\delta$ unital and   & $J=(\ker\delta)^\bot$                 & covariance algebra \cite{kwa}
			\\ & $A$ commutative        &                                       &
			\\ \hline			5. & arbitrary  & $ \{0\}\subset  J \subset (\ker\delta)^\bot$       & partial-isometric 
			\\ &            &                                                    & crossed product	\cite{kwa-leb},
									\\ \hline
			6. & arbitrary  & $J=A$                 & isometric crossed product		
			\\ &            &                       & \cite{Adji_Laca_Nilsen_Raeburn}				
			\\ \hline
			7. & arbitrary  & $J=\{0\}$                 & partial-isometric 		
			\\ &            &                           & crossed product	\cite{Lin-Rae}
			\\ \hline
		\end{tabular}
		\caption{Different crossed products as relative Cuntz-Pimsner algebras}
\end{center}
\end{table}
\par
%Let us start with showing that  a reduction of a Hilbert bimodule associated with a $C^*$-dynamical system leads again to a  Hilbert bimodule of another ("smaller") $C^*$-dynamical system.
%\\
The first step is to show that the reduction procedure presented in  Section \ref{Reduction of Hilbert bimodule} when started with $C^*$-dynamical system bimodule leads again to another $C^*$-dynamical system bimodule.
Indeed, let $E$ be a Hilbert bimodule of the $C^*$-dynamical system $(\A,\delta)$ and let $J$ be an ideal in $A$ and  $J_n$, $n=0,1,...,\infty$, the related ideals  defined in Section \ref{Reduction of Hilbert bimodule}. Then 
$$
J_n=\underbrace{\delta^{-1}\Big(\delta^{-1}\big(...(\delta^{-1}}_{n\textrm{ times}}(\ker\delta \cap J)\cap J)...\big)\cap J\Big)\cap J , \,\,\,\qquad n\in \N.
$$
That is 
$$
J_n=\delta^{-n}(\ker \delta)\cap \bigcap_{k=0}^n \delta^{-k}(J)
$$
and hence 
\begin{equation}\label{ideal I_infinity}
J_\infty=\overline{\{a \in J: \exists_{n\in\N}\,\, \delta^n(a)=0)\}} \cap \bigcap_{n\in \N} \delta^{-n}(J).
\end{equation}
Let $n\in \N\cup\{\infty\}$. Since the ideal $J_n$ is  $E$-invariant,   taking $x=\delta(1)$ in \eqref{left action} one sees  that the    mapping $\delta_n:A/J_{n} \to A/J_n$  given by 
$$
\delta_n\circ q^{J_n} =q^{J_n}\circ \delta 
$$
is a well  defined endomorphism of $A/J_n$ and the quotient Hilbert bimodule $E/E J_n$ may be viewed as the $C^*$-dynamical system Hilbert bimodule of $(A/J_n, \delta_n)$. In particular
\begin{equation}\label{delta_infinity}
\delta_\infty(a+ J_{\infty}):=\delta(a) + J_\infty
\end{equation}
is an endomorphism of $A/J_{\infty}$ such that $\ker \delta_\infty \cap q^{J_{\infty}}(J)=\{0\}$.
\\
  Obviously one may apply  Theorem \ref{reduction thm} to  each of the systems $(A/J_n,\delta_n)$, $n\in \N\cup\{\infty\}$, however, we  focus on the case $n=\infty$. Then by virtue of Proposition \ref{universality proposition} we get
\begin{prop}\label{reducing C*-Hilbert bimodules}
If  $E$ is a Hilbert bimodule of the $C^*$-dynamical system $(A,\delta)$ and  $J$ is an ideal in $A$, then 
$$\OO(J,E)=\OO( q^{J_\infty}(J), E/EJ_\infty)$$
is  a universal algebra generated by a copy of the algebra $A/J_{\infty}$ and a partial isometry $u$ subject to relations 
$$
ua u^*=\delta_\infty(a),\,\,\, a\in A/J_{\infty}, \qquad u^*u\in (A/J_{\infty})',
$$
$$
q^{J_\infty}(J)=\{a\in A/J_{\infty}: u^*u a=a\}.
$$
\end{prop}

\begin{cor}\label{kernel of a covariant representation}
If $(\pi, U)$ is a covariant representation of $(\A,\delta)$ associated with an ideal $J\in A$, then 
$
J_\infty \subset \ker \pi. 
$ 
\end{cor}

\begin{cor}\label{corl for reduced dynamics}
The relative Cuntz-Pimsner algebra $\OO(J,E)$  coincides with  crossed product N.6, Table 1, applied to the $C^*$-dynamical system $(A/J_\infty,\delta_\infty)$, and the ideal $q^{J_\infty}(J)$.
\end{cor}
\begin{cor}
The Cuntz-Pimsner algebra $\OO(A,E)$ (crossed product N.6, Table 1) reduces to zero if and only if the set of elements $a\in A$ such that $\delta^n(a)=0$ for certain $n\in\N$, is dense in $A$.
\end{cor}
Above results show us how to  reduce the investigation of crossed products to the case where $J$ satisfies  \eqref{inclusion for an ideal J}, thereby let us assume for a while that  $(A,\delta)$ be a $C^*$-dynamical system and  $J$ is an ideal in $A$ such that $\ker\delta\cap J=\{0\}$. As in \cite{kwa-leb} we slightly extend  $(A,\delta)$ to a certain system $(A_J,\delta_J)$ for which we shall have $\OO(J,E)=\OO((\ker\delta_J)^\bot,E_J)$. For this purpose we denote by $A_J$ the direct sum of quotient algebras
$$
A_J=\big(A/\ker\delta \big) \oplus \big(A/J\big),
$$
and  we set $\delta_J:A_J\to A_J$ by the formula
\begin{equation}
\label{d_J}
 A_J\ni \big((a +\ker\delta)\oplus (b
+J)\big)\stackrel{\delta_J}{\longrightarrow}(\delta(a) +\ker\delta)\oplus
(\delta(a) +J)\in A_J.
\end{equation}
Endomorphism  $\delta_J$ is  well defined and since $
\ker \delta_J = 0\oplus A/J 
$ its kernel is unital. Moreover,  $A$ embeds  into $C^*$-algebra $A_J$ via
\begin{equation}
\label{A_in} A \ni a \longmapsto \big(a +\ker\delta\big)\oplus \big(a +
J\big)\in A_J.
\end{equation}
 Since $\ker\delta \cap J = \{0\}$ this mapping is injective  and we may
treat $A$ as the corresponding subalgebra of $A_J$. Under this
identification $\delta_J$ is an extension of $\delta$.
\begin{defn}\label{definition of the canon}
Let $(A,\delta)$ be a $C^*$-dynamical system and  $J$  an arbitrary ideal in $A$. Let $((A/J_\infty)_{q^{J_\infty}(J)},(\delta_\infty)_{q^{J_\infty}(J)})$ be the above constructed extension of the reduced $C^*$-dynamical system $(A/J_\infty,\delta_\infty)$  given by \eqref{ideal I_infinity}, \eqref{delta_infinity}. We shall write 
$$
(A_J,\delta_J):=(A/J_\infty)_{q^{J_\infty}(J)},(\delta_\infty)_{q^{J_\infty}(J)})
$$
and say that $(A_J,\delta_J)$ is the \emph{canonical $C^*$-dynamical system} associated to $(A,\delta)$ and $J$.
\end{defn}
Combining Proposition \ref{reducing C*-Hilbert bimodules} with \cite[Proposition 1.2]{kwa-leb}, see also \cite[Corollary 1.7]{kwa-ext}, we get the following 
\begin{thm}
Let   $E$ be a Hilbert bimodule of the $C^*$-dynamical system $(A,\delta)$ and let $J$ be an ideal in $A$. If $E_J$ is the  Hilbert bimodule of the canonical system $(A_J,\delta_J)$, then 
$$\OO(J,E)=\OO( (\ker\delta_J)^\bot, E_J)$$
is  a universal algebra generated by a copy of the algebra $A_J$ and a partial isometry $u$ subject to relations 
\begin{equation}\label{reduced realtions}
u a u^*=\delta_J(a),\,\,\, a\in A_J, \qquad u^*u\in A_J
\end{equation}
(relations \eqref{reduced realtions} imply that $u^*u$ belongs to the center of $A_J$).
\end{thm}

The usefulness of canonical $C^*$-dynamical system $(A_J,\delta_J)$ manifests in reducing fairly  complicated relations  \eqref{covariance rel1}, \eqref{covariance rel2}, \eqref{covariance rel3} which may degenerate to the nondegenerate relations \eqref{reduced realtions}.  Moreover, one may apply  the results of \cite{kwa-leb} and \cite{kwa} to  $(A_J,\delta_J)$ (in particular, norm evaluation of elements \cite[Section 3]{kwa-leb}, isomorphism theorem \cite[Section 5]{kwa-leb}, \cite[Subsection 6.2]{kwa}, ideal structure \cite[Subsection 6.1]{kwa}) and thus get   results concerning relative Cuntz-Pimsner algebras $\OO(J,E)$ of $C^*$-dynamical systems Hilbert bimodules. 
\\
In fact, one could go even further and use the construction  from \cite{kwa-ext} to  extend,  the canonical system $(A_J,\delta_J)$ up to a $C^*$-dynamical system $(B,\tdelta)$ possessing  a complete transfer operator. Then $B$  corresponds to the fixed point subalgebra of $\OO(J,E)$ for the gauge action $\gamma$  (Proposition \ref{RCP algebra}),  and one could  apply  the results of \cite{Ant-Bakht-Leb} or  isomorphism theorem \cite[Section 6]{kwa-int} to $(B,\tdelta)$ to study $\OO(J,E)$ in terms of 'Fourier' coefficients.

\noindent \textsc{Bartosz K. Kwa\'sniewski}\\
Institute of Mathematics,  University  of Bialystok,\\
 ul. Akademicka 2, PL-15-267  Bialystok, Poland \\
 \emph{e-mail:} \texttt{bartoszk@math.uwb.edu.pl}\\
 \emph{www:} \texttt{http://math.uwb.edu.pl/$\sim$zaf/kwasniewski}
 \bigskip

\noindent \textsc{Andrei V. Lebedev}\\
Institute of Mathematics,  University  of Bialystok,\\
 ul. Akademicka 2, PL-15-267  Bialystok, Poland\\
 \emph{e-mail:} \texttt{lebedev@bsu.by}

\end{document}